\documentclass[12pt]{amsart}
\usepackage[all]{xy}
\usepackage{amssymb}
\usepackage{amsthm}
\usepackage{amsmath}
\usepackage{amscd,enumitem}
\usepackage{verbatim}
\usepackage{eurosym}
\usepackage{graphicx}
\usepackage{dsfont}
\usepackage{stmaryrd}
\usepackage{color}
\usepackage{float}
\usepackage{longtable}
\usepackage{dcolumn}
\usepackage[mathscr]{eucal}
\usepackage[all]{xy}
\usepackage{hyperref}
\usepackage[title]{appendix}
\usepackage[dvipsnames]{xcolor}
\usepackage{bbm}
\usepackage[textheight=8.85in, textwidth=6.8in]{geometry}
\usepackage{multirow}
\usepackage{caption}
\theoremstyle{plain}
\newtheorem{theorem}{Theorem}[section]
\newtheorem{corollary}[theorem]{Corollary}
\newtheorem{lemma}[theorem]{Lemma}
\newtheorem{proposition}[theorem]{Proposition}

\theoremstyle{definition}
\newtheorem{definition}[theorem]{Definition}

\theoremstyle{remark}

\newcommand{\Q}{\mathbb{Q}}

\newcommand{\F}{\mathbb{F}}

\newcommand{\C}{\mathbb{C}}

\newcommand{\leg}[2]{\genfrac{(}{)}{}{}{#1}{#2}}

\catcode`,\active

\catcode`\,12

\numberwithin{equation}{section}

\begin{document}
	\title[]{Higher moments of twisted Kloosterman sums and the arithmetic of algebraic varieties}
		\author{Manash Dutta and Neelam Saikia}

	\address{School of Basic Sciences, Indian Institute of Technology Bhubaneswar, Argul, Khordha 752050, Odisha, India.}
	\email{a25ma09002@iitbbs.ac.in}

	\address{School of Basic Sciences, Indian Institute of Technology Bhubaneswar, Argul, Khordha 752050, Odisha, India.}
	\email{neelamsaikia@iitbbs.ac.in}

	
	\begin{abstract}
     For an odd prime $p,$ consider the affine varieties $X_k: \sum\limits_{i=1}^{k}(x_i+x_i^{-1})=0$ over the finite field $\F_p.$ In this paper, a key feature of our work is arithmetic interpretation between the higher power moments of Kloosterman sums, modular forms and arithmetic geometry of the varieties $X_5$ and $X_6$ through their rational point counts. Our main results provide two distinct explicit formulas for the point count $\#X_5(\F_p).$ In the first formula, we express $\#X_5(\F_p)$ in terms of Fourier coefficients of modular forms and a certain weighted average of  fourth power twisted Kloosterman moments, where the average is taken over multiplicative characters. We further establish an identity expressing this weighted average of fourth power moment of twisted Kloosterman sums in terms of fifth power quadratic twisted Kloosterman moment, together with Fourier coefficients of modular forms. Moreover, we derive an independent expression for the same  average with different weight directly in terms of Fourier coefficients of modular forms. In the second identity for $\#X_5(\F_p)$, we express this purely in terms of Fourier coefficients of two modular forms. Moreover, we express $\#X_6(\F_p),$ in terms of Fourier coefficients of three modular forms. For these two expressions of $\#X_5(\F_p)$ and $\#X_6(\F_p),$ our approach is based on an identity that relates the point counts of $X_k$ to the $k$-th power moments of Kloosterman sum and its quadratic twist.

         \end{abstract}

	\maketitle
	\section{Introduction and Statement of results}

Let $p$ be an odd prime and $\F_p$ denote a finite field with $p$ elements. We write $\zeta_p:=e^{2\pi i/p}$ and for $a\in\F_p$ consider the Kloosterman sum
\begin{equation}\label{K(a)}
K(a):=\sum\limits_{x\in\F_p^{\times}}\zeta_p^{x+\frac{a}{x}}.
\end{equation}
For $n\geq1$, let $S(n)$ denote the power moment
\begin{equation}\label{132}
    S(n):=\sum\limits_{a=1}^{p-1}K(a)^n .
\end{equation}
Following \cite{iwaniec}, we state certain known values for $S(n)$ as given by
\begin{equation}\label{S(3)}
    S(1)=1, \ S(2)=p^2-p-1, \  S(3)=\leg{p}{3}p^2+2p+1,
\end{equation}

where $\leg{\cdot}{p}$ denotes the Legendre symbol modulo $p,$ and  
\begin{equation}\nonumber
    S(4)=2p^3-3p^2-1.
\end{equation}
 Following \cite{Evans, Livne, Peters},for $p>5,$ we write
\begin{equation}\label{s(5)}
    S(5)=\leg{p}{3}4p^3+(a_g(p)+5)p^2+4p+1,
\end{equation}
where $a_g(p)$'s are the $p$-th Fourier coefficient of the $q$-expansion of the newform, 
\begin{equation}\label{g(z)}
g(z):=\eta(z)\eta(3z)\eta(5z)\eta(15z)(\theta_2(q)\theta_2(q^{15})+\theta_3(q)\theta_3(q^{15})),
\end{equation} 
of weight $3$ on $\Gamma_0(15)$ with quadratic nebentypus of conductor $15,$ where 
$\eta$ is the Dedekind eta function, $\theta_2,\theta_3$ are the classical Jacobi theta functions and $q=e^{2\pi i z}.$ Moreover, for $p>5,$ the Fourier coefficients $a_g(p)$'s are given by
\begin{equation}\label{a_p}
    a_g(p)=\begin{cases}
        2p-12u^2,& \ \text{if}\ p=3u^2+5v^2;\\
        4x^2-2p,&\ \text{if}\ p=x^2+15y^2;\\
        0,&\ \text{if}\ p\equiv7,11,13,\ \text{or}\ 14 \pmod{15}.
    \end{cases}
\end{equation}

Furthermore, if $p>7,$ then due to the work of \cite{Hulek}, and \cite[Eq. (1.9)]{Evans}, it follows that
\begin{equation}\label{S(6)}
    S(6)=5p^4-10p^3-(a_u(p)+9)p^2-5p-1,
\end{equation}
where $a_u(p)$'s are $p$-th the coefficient of the $q$-expansion of the newform of weight 4 and level 6 given by  \begin{equation}\label{u(z)}
    u(z):=(\eta(6z)\eta(3z)\eta(2z)\eta(z))^2,
\end{equation} where $\eta$ is the Dedekind eta function.
To this end, we recall the twisted power moments of Kloosterman sums. If $\chi$ is a multiplicative character of $\F_p$ with the convention that $\chi(0)=0$ for all multiplicative characters $\chi$ including the trivial character $\varepsilon,$ then 
$n$-th power moment of twisted Kloosterman sum is denoted by $S(n,\chi)$ and defined by 
\begin{equation}\label{Twisted-Moments}
    S(n,\chi):=\sum_{a=1}^{p-1}\chi(a)K(a)^n.
\end{equation}
In particular, if we put $n=1$ in \eqref{Twisted-Moments}, then we can easily verify that $S(1,\chi)=g(\chi)^2,$ where $g(\chi)$ is Gauss sum defined by
\begin{equation}\label{Gauss-sum}
    g(\chi):=\sum_{x=0}^{p-1}\chi(x)\zeta_p^x.
\end{equation} Moreover, if $n=2$ and $\chi\neq\varepsilon,$ then Liu \cite{Liu} showed that 
\begin{equation}\label{eqn-Liu}
    S(2,\chi)=p\phi(-1)g(\chi)^2{\chi}(1/4)\binom{\chi}{\phi},
\end{equation}
where $\phi$ is quadratic character of $\F_p$ and $\binom{\chi}{\phi}$ is the normalized Jacobi sum defined in \eqref{Jacobi-1}. For a multiplicative character $\chi,$ $n$-th twisted Kloosterman sheaf sum of $\chi$ is defined by 
\begin{equation}\label{T(n)}
   T(n,\chi):=\sum\limits_{a\in\F_p^{\times}}\chi(a)(\alpha(a)^n+\alpha(a)^{n-1}\beta(a)+\alpha(a)^{n-2}\beta(a)^2+\cdots+\alpha(a)\beta(a)^{n-1}+\beta(a)^{n}), 
\end{equation}
where $\alpha(a),\ \beta(a)$ satisfy $x^2+K(a)x+p=0.$ Thus,  $\alpha(a)+\beta(a)=-K(a)$ and  $\alpha(a)\cdot\beta(a)=p.$ Kloosterman sums have
significant consequences on the theory of automorphic forms, Galois representations and related areas, see
for example \cite{Evans,Davenport1,Deshouillers,Evan1,Hulek,iwaniec,Livne,Moisio,Salie,Schmidt}. The interplay between the arithmetic of Calabi-Yau varieties, modular forms, and Kloosterman sums provides a powerful source of explicit arithmetical formulas. For example, Ahlgren and Ono \cite{ahlgren-ono} discussed the modularity of the Calabi-Yau threefold 
$$ X_4:x_1+x_1^{-1}+x_2+x_2^{-1}+x_3+x_3^{-1}+x_4+x_4^{-1}=0 ,$$ by establishing the number of $\F_p$-rational points of $X_4$ in terms of Fourier coefficient of weight 4 newform. In this direction, the second author \cite{Saikia} related the number of $\F_p$-rational points of the variety $X_4$ to the fourth power moment of twisted Kloosterman sum. In the spirit of the work in \cite{ahlgren-ono}, our main purpose of this paper is to establish arithmetic connections between higher moments of Kloosterman sums, modular forms, and the point counts of higher dimensional Calabi-Yau type varieties. 
More precisely, we consider the following $(k-1)$-dimensional variety defined by
$$
X_k:\sum\limits_{i=1}^k (x_i+x_i^{-1})=0.
$$
If we set $k=4,$ then $X_4$ is the Calabi-Yau threefold studied in \cite{ahlgren-ono, Saikia}. Therefore, it is natural to seek similar questions if $k>4.$
To this end, we define 
$$
 N_k(p):=\#X_k(\F_p)=\#\{(x_1,x_2,\ldots,x_k)\in({\F}_p^{\times})^k:
\sum\limits_{i=1}^{k}(x_i+x_i^{-1})=0\}. 
$$
 As our first theorem, we express $N_5(p)$ in terms of Fourier coefficients of two modular forms and a weighted average of fourth moment of twisted Kloosterman sums with weights given by the reciprocal of the squares of the Gauss sums associated with the quadratic twisted character $\phi\chi$, where the average is taken over multiplicative characters $\chi.$ To state this result, we fix some notation. For $n\geq1,$ suppose that $a_h(n)$ are the $n$-th Fourier coefficients of the $q$-expansion of $h(z),$ a newform of weight $4$ and level $8,$ defined by
\begin{equation}\label{b_p}
h(z):=\sum\limits_{n=1}^{\infty}a_h(n)q^n=q\prod\limits_{n=1}^{\infty}(1-q^{2n})^4(1-q^{4n})^4,
\end{equation}
where $q=e^{2\pi iz}.$
\begin{theorem}\label{theorem-1}
    For a prime $p>5,$ we have  \begin{align*} N_5(p)=&\frac{(p-1)^5+1}{p}-\frac{2p(1+\phi(-1))}{p-1}+pa_g(p)+4p^2\phi(-3)\\
    &+\frac{(p+2)(2p-1)}{p-1}+\frac{p}{p-1}\sum\limits_{\psi\neq\varepsilon,\phi}\frac{S(4,\psi)}{g(\phi\psi)^2}+T,
   \end{align*}
   where $$  T:=\begin{cases}
-2p^2+5p+1-\frac{a_h(p)+p^2}{p-1}, &\text{if} \ p\equiv 3 \pmod{4};\\
\frac{2p^3-6p-3-a_h(p)}{(p-1)}, & \text{if} \ p\equiv 1 \pmod{4},
\end{cases}$$
and where $a_g(p),\ \text{and}\ a_h(p)$ denote the $p$-th Fourier coefficients of the modular forms $g$ and $h,$ of weights $3$ and $4,$ respectively, as defined in \eqref{a_p}, and \eqref{b_p}. 
\end{theorem}

In the next theorem, we express the weighted average of   fourth twisted Kloosterman moments that arise in the expression of $N_5(p)$ as given in the above theorem in terms of fifth power moment of quadratic twisted Kloosterman sum and $p$-th Fourier coefficients of modular form $h,$ of weight $4$ and level $8$ as defined in \eqref{b_p}.

 \begin{theorem}\label{corollary-1} For prime $p>5$, we have
    \begin{align*}
   \frac{p^2}{p-1}\sum\limits_{\psi\neq\varepsilon,\phi}\frac{S(4,\psi)}{g(\phi\psi)^2}=S(5,\phi)+\frac{p(3p^2-4p-2)}{p-1}+\frac{2p^2(1+\phi(-1))}{p-1}-p\cdot T,
    \end{align*}
    where  $T$ is same as defined in Theorem \ref{theorem-1} and $a_h(p)$ that arise in the expression of $T$ is the $p$-th Fourier coefficient of $h$ as given in \eqref{b_p}.
    \end{theorem}

In the following theorem, we further express the weighted average of fourth twisted Kloosterman moments, with weights given by the reciprocal squares of the corresponding Gauss sums in terms of $p$-th Fourier coefficients of modular forms $g$ and $h.$

\begin{theorem}\label{theorem-5}
    For prime $p>5,$ we have 
$$\sum\limits_{\psi\neq \varepsilon,\phi}\frac{S(4,\psi)}{g(\psi)^2}=(p-1)\left(4p\phi(-3)+\frac{p+2}{p-1}+\frac{2p(1+\phi(-1))}{p-1}+a_g(p)+\frac{\phi(-1)}{p-1}(a_h(p)+p)\right),$$
where $a_g(p)$ and $a_h(p)$ are the $p$-th Fourier coefficients of the modular forms $g$ and $h,$ respectively (see \eqref{a_p} and \eqref{b_p}).
\end{theorem}

We now turn to another main objective of the paper. In the following results, we derive an alternative expression for $N_5(p)$ purely in terms of $p$-th Fourier coefficients of two modular forms, as well a new formula for $N_6(p)$ involving $p$-th Fourier coefficients of three modular forms. Furthermore, we establish a modular form interpretation of the weighted average of the fourth twisted Kloosterman moment appearing in the first expression for $N_5(p),$ as given in Theorem \ref{theorem-1}. To state these results more precisely, we fix some notation. Let $S_k(\Gamma_0(N))$ denote the space of cusp forms of weight $k$ and level $N$. Suppose that $f$ and ${\tilde{f}}$ are two modular forms with the following $q$-expansions:
    $$
f(z)=\sum_{n=1}^{\infty}a_f(n)q^n \ \text{and} \ \ {\tilde{f}}(z)=\sum_{n=1}^{\infty}a_{\tilde{f}}(n)q^n,
$$
where $q=e^{2\pi iz}$ and with the convention that throughout these two modular forms are normalized cuspidal Hecke eigenforms. Thus, we have $a_f(1)=a_{\tilde{f}}(1)=1$.
      Then we fix $f$ and ${\tilde{f}}$ as the following forms, identified by their labels in the LMFDB as referred to the online database \cite{LMFDB} of $L$-functions and modular forms by setting
\begin{align}
&f:=f_{60.3.b.a}\in S_3(\Gamma_0(60)),\ \text{and}\label{121}\\
&\tilde{f}:={\tilde{f}}_{24.6.a.a}\in S_6(\Gamma_0(24))\label{122}.
\end{align}
 For $f_{60.3.b.a}$, the Fourier coefficients lie in $\mathbb{Q}(\sqrt{5},i).$ 
The arithmetic of Calabi–Yau varieties over finite fields has interesting connections with the theory of modular forms. In particular, the number of rational points over finite fields can sometimes be expressed in terms of Fourier coefficients of modular forms. For example, see \cite{ahlgren-ono}. In the next two theorems, we present two such connections for $X_5$ and $X_6$. More precisely, in Theorem \ref{theorem-3} we relate $N_5(p)$ and the $p$-th Fourier coefficients of the modular forms $f$ and $g$ as defined in \eqref{121} and \eqref{g(z)}, respectively. In Theorem \ref{theorem-4}, we express $N_6(p)$ in terms of the $p$-th Fourier coefficients of the modular forms $u$, $h,$ and $\tilde{f}$ as defined in \eqref{u(z)}, \eqref{b_p}, and \eqref{122}, respectively.

 \begin{theorem}\label{theorem-3}
 Let $p>5$ be a prime. Then we have 
  \begin{align*} N_5(p)=&p^4-5p^3+10p^2+9+(p^2+a_f(p)^2)\cdot\phi(-1)\\
    &+4p^2\phi(-3)+(a_g(p)-5)p-\phi(15)p^2+4p^3\cdot c_p ,
\end{align*} 
where $a_f(p),$ and $a_g(p)$ are the $p$-th Fourier coefficients of $f,$ and $g$ as defined in \eqref{121}, and \eqref{a_p}, respectively. Moreover, \begin{equation}\label{c_p}
    c_p:=\begin{cases}
    -\frac{\phi(3)}{p}, \ &\text{if}\  p=2\pmod3;\\
    \frac{\phi(3)(4x^2-p)}{p^2},\ &\text{if}\ p=1\pmod3 \text{and}\ x^2+3y^2=p.
\end{cases}
\end{equation}
 \end{theorem}
 
\begin{theorem}\label{theorem-4}
 Let $p>7$ be a prime and then we have the following:
 $$N_6(p)=p^5-6p^4+20p^3-40p^2-p(5\cdot a_h(p)+a_u(p)-6)-11-a_{\tilde{f}}(p),$$
where $a_h(p),a_u(p),$ and $a_{\tilde{f}}(p)$ are the $p$-th Fourier coefficients of $h,$  $u,$ and $\tilde{f}$ as defined in \eqref{b_p}, \eqref{u(z)}, and \eqref{122}, respectively.
\end{theorem}
Now, combining Theorem \ref{theorem-1}
 and Theorem \ref{theorem-3}, we obtain the following result expressing the weighted average of fourth power moment of twisted Kloosterman sums in terms of the Fourier coefficients of modular forms $f$ and $h.$
 \begin{corollary}\label{corollary-2} For prime $p>5$, we have
   \begin{align*}
     \frac{p}{p-1}\sum\limits_{\psi\neq\varepsilon,\phi}\frac{S(4,\psi)}{g(\phi\psi)^2}=& 4p^3\cdot c_p+(p^2+a_f(p)^2)\phi(-1)-\phi(15)p^2 \\&+\frac{2p(1+\phi(-1))}{p-1}+\frac{3p^2-4p-2}{p-1}-T
    ,\notag
\end{align*}
where $a_f(p),$ and $a_h(p)$ are the $p$-th Fourier coefficients of $f,$ and $h$ as defined in \eqref{121} and \eqref{b_p}, respectively. Moreover $c_p,$ and $T$ are same as defined in \eqref{c_p}, and Theorem \eqref{theorem-1}, respectively.
 \end{corollary}

The rest of the paper is structured as follows: In Section 2, we recall some elementary results for multiplicative characters, Gauss and Jacobi sums. Furthermore, in Section 2, we introduce Gaussian hypergeometric functions over finite fields introduced by Greene \cite{Greene}. This function plays an important role in the proofs of our main results. In Section 3, we prove certain propositions that are key ingredients to prove our main results. Finally, we provide the proofs of Theorem \ref{theorem-1} - Theorem \ref{theorem-4} in Section 4.

\section{Notation and Preliminary results}
In this section, we recall certain basic definitions including Jacobi sums and elementary results. We begin with multiplicative characters.
\subsubsection*{Multiplicative characters:} A multiplicative character $\chi$ on $\F_p^\times$ is a group homomorphism $\chi:{\F}_p^\times\rightarrow {\C^{*}}$, which can be extended to ${\F}_p$ by simply setting $\chi(0):=0$ for all multiplicative characters, including the trivial character $\varepsilon.$ We denote the group of multiplicative characters on ${\F}_p $ by $\widehat{{\F}_p^{\times}}.$
For two multiplicative character $A,B\in\widehat{\mathbb{F}_p^{\times}}$, Jacobi sum is defined as 
\begin{equation}\label{jacobi-3}
    J(A,B):=\sum_{x=0}^{p-1}A(x)B(1-x),
\end{equation}
and the binomial $\binom{A}{B}$ is normalized Jacobi sum as defined by
\begin{equation}\label{Jacobi-1}
   \binom{A}{B}:=\frac{B(-1)}{p}J(A,\overline{B}) ,
\end{equation}
where $\overline{B}$ is the inverse of $B.$
The following lemma gives orthogonality results for multiplicative characters.

\begin{lemma}\cite[Chapter $8$]{rosen}
    For a character $\chi\in\widehat{{\F}_p^{\times}},$ we have  
\begin{equation}\label{ortho-1}
\sum\limits_{x=0}^{p-1}\chi(x)=\begin{cases}
    p-1, &  \text{ if}\  \chi=\varepsilon ;\\ 
    0,&  \text{ if}\  \chi\neq\varepsilon ,\\
\end{cases}
\end{equation}
and \begin{equation}\label{ortho-2}
    \sum\limits_{\chi\in\widehat{{\F}_p^{\times}}}\chi(x)=\begin{cases}
    p-1,&  \text{ if}\ x=1;\\
    0,&  \text{ if}\ x\neq1.
\end{cases}
\end{equation}
\end{lemma}

\begin{lemma}\cite[Lemma 2.1]{ahlgren-ono}\label{lem-2.1}
Let $p$ be an odd prime and suppose that $a\in\F_p$. Then we have,
$$\sum\limits_{x\in\F_p}\phi(x^2+ax)=\sum\limits_{x\in\F_p}\phi(x^2+a)=\begin{cases}
    p-1, & \ \text{if}\ a=0, \\
    -1,  & \ \text{if}\ a\neq0.
\end{cases}$$
    
\end{lemma}
Now, we recall certain elementary identities of Gauss and Jacobi sums from \cite{rosen}.
\begin{equation}\label{Gauss}
    g(A)g(\overline{A})=pA(-1)-(p-1)\delta(A),
\end{equation}
and \begin{equation}\label{Jacobi-2}
    J(A,B)=\frac{g(A)g(B)}{g(AB)}+(p-1)B(-1)\delta(AB),
\end{equation}
where $\delta(\chi)=\begin{cases}
    0,& \ \text{if}\ \chi\neq\varepsilon;\\
    1,&\ \text{if}\ \chi = \varepsilon.
\end{cases}$

From \cite{Davenport}, we state the following relation as a special case of Hasse-Davenport identity:
\begin{equation}\label{Hasse}
g(A^2)g(\phi)\overline{A}(4)=g(A)g(\phi A).
\end{equation}
We now state the following lemma that gives some cancellations:
\begin{lemma}\label{lem-3} 
For $A\neq\varepsilon,B\neq\varepsilon,\text{and}\ A\overline{B}\neq\varepsilon$  then,
$$\binom{A}{B}\binom{B}{A}=\frac{AB(-1)}{p}.$$
\end{lemma}

\begin{proof}
  The proof follows directly by using \eqref{Jacobi-1}, \eqref{Gauss} and \eqref{Jacobi-2}.
\end{proof}

\noindent We now recall $_2F_1$-Gaussian hypergeometric functions introduced by Greene \cite{Greene} over a finite field $\F_p.$
\begin{definition}\label{hypergeometric-2}\cite[Definition $3.5$]{Greene} For characters $A,B$ and $C$ of $\F_p$ and $x\in\F_p$,
   $$ {_2 F_1}\left(\begin{matrix}
        A,&B\\
        &C
    \end{matrix}
    \middle|x
    \right):=\varepsilon(x)\frac{BC(-1)}{p}\sum\limits_{y\in\F_p}B(y)\overline{B}C(1-y)\overline{A}(1-xy).
$$
\end{definition}
This function can be described as a finite field analogue of classical hypergeometric series studied by Gauss, Euler, Ramanujan etc.
Moreover, Greene \cite{Greene} extended this function to $_{n+1}F_n$-Gaussian hypergeometric functions for any $n\geq1$ over finite fields involving normalized Jacobi sums as follows:
\begin{definition}\cite[Definition $3.10$]{Greene}\label{hypergeometric-1} For multiplicative characters $A_0,A_1,..,A_n,B_1,B_2,...,B_n$ of ${\F}_p$ and $x\in {\F}_p,$  $_{n+1}F_n$-Gaussian hypergeometric function is defined by
$$_{n+1}F_n\left(\begin{matrix}
A_0,&A_1,&\cdots,&A_n\\
& B_1,&\cdots,&B_n\\
\end{matrix}
\middle|x
\right):=\frac{p}{p-1}\sum\limits_{\chi\in\widehat{{\F}_p^{\times}}}\binom{A_0\chi}{\chi}\binom{A_1\chi}{B_1\chi}\cdots\binom{A_n\chi}{B_n\chi}\chi(x).$$
\label{def-1}\end{definition}
Many authors have established that these functions exhibit several properties analogous to those of the classical hypergeometric series studied by Gauss and Kummer. They have found important applications in the study of Ap\'{e}ry-style supercongruences, the Eichler–Selberg trace formula, Galois representations, and zeta functions of arithmetic varieties. Another fundamental direction is investigating special values of these hypergeometric functions. For example, Greene \cite[Theorem 4.9]{Greene} derive that \begin{equation}\label{2F1}{_2F_1}\left(\begin{matrix}
        A,&B\\
        &C
    \end{matrix}
    \middle|1
    \right)=A(-1)\binom{B}{\overline{A}C}.\end{equation}
    Furthermore, in this direction, for prime $p>3$, Ono \cite[Theorem 6 (v)]{ono} showed 
\begin{align}\label{3F2}
{_3 F_2}\left(\begin{matrix}
        \phi,&\phi, & \phi\\
        &\varepsilon, & \varepsilon
    \end{matrix}
    \middle|\frac{1}{4}
    \right)=\begin{cases}
    -\frac{\phi(3)}{p}, \ &\text{if}\  p=2\pmod3;\\
    \frac{\phi(3)(4x^2-p)}{p^2},\ &\text{if}\ p=1\pmod3\ \text{and}\ x^2+3y^2=p.
\end{cases}
\end{align}
In the next section, we heavily use Gaussian hypergeometric functions as our tools to prove the main results.
 
\section{Arithmetic of Gaussian hypergeometric functions}
In this section, we derive some important propositions that are required to prove the main results. Our main goal of this section is to write $N_5(p)$ as a character sum involving $_4F_3$-hypergeometric functions that is defined in the previous section. To do this we need to derive series of propositions and lemmas. We begin with the following proposition:
\begin{proposition}\label{proposition-1}
For a prime $p>3$ we have $N_5(p)=S+p^4-5p^3+10p^2-10p+5$,
where $$S:=\sum\limits_{c,d}\phi((c+d)^2-1)F(c)F(d) \ \text{and}\  F(t):=\sum\limits_{x}\phi(x^2-1)\phi((x+t)^2-1).$$
\end{proposition}

\begin{proof} We may write \begin{align}
N_5(p)&=  \sum\limits_{a,b,c,d\in\F_p} \#\ \{x_1+x_1^{-1}=a \}\  \#\{x_2+x_2^{-1}=b \}\ \#\{x_3+x_3^{-1}=c \}\ \#\left\{x_4+x_4^{-1}=d \right\} \notag\\
&\times \# \left\{x_5+x_5^{-1}=-a-b-c-d \right\}.\notag
\end{align}
Observe that $x_1+{x_1^{-1}}=a$ if and only if $\left(x_1-\frac{a}{2}\right)^2=\frac{a^2}{4}-1.$ Now, using this fact and replacing $a$ by $2a,$ we have $$\#\{x_1+x_1^{-1}=a \}=\phi(a^2-1)+1.$$ This yields,
\begin{align*}
N_5(p)=&\sum\limits_{a,b,c,d\in\F_p}{\{ \phi(a^2-1)+1\}\{ \phi(b^2-1)+1\}\{ \phi(c^2-1)+1\}\{ \phi(d^2-1)+1\}}\\
&\times \{ \phi((a+b+c+d)^2-1)+1\}.
\end{align*}

Now, in the above expression of $N_5(p),$ if we replace $c$ by $c-b,$ $d$ by $d-a,$ and subsequently $b$ by $-b,$ and $a$ by $-a,$ then after rearranging the terms, we have
\begin{align*}
N_5(p) &=  \sum\limits_{a,b,c,d\in\F_p} {\{\{ \phi(a^2-1)+1\}\{\phi((d+a)^2-1)+1 \} \}\{\{ \phi(b^2-1)+1\}\{\phi((b+c)^2-1)+1 \} \}}\\
 & \times \{\phi((c+d)^2-1)+1 \} \\
&=  \sum\limits_{d}\left(F(d)+\sum\limits_{a}\phi(a^2-1)+\sum\limits_{a}\phi((a+d)^2-1)+p \right)\\  
& \times \sum\limits_{c} \left(F(c)+\sum\limits_{b}\phi(b^2-1)+\sum\limits_{b}
\phi((b+c)^2-1)+p \right)  \{\phi((c+d)^2-1)+1 \}.
\end{align*}
By using Lemma \ref{lem-2.1}, we obtain that
\begin{align*}
\sum\limits_{d}F(d)&=\sum\limits_{d}\sum\limits_{a}\phi(a^2-1)\phi((a+d)^2-1)= \sum\limits_{a}\phi(a^2-1)\sum\limits_{d}\phi((a+d)^2-1)=1.
\end{align*}
Similarly, we have
\begin{equation}\label{2.18}
    \sum\limits_{c}F(c)=1=\sum\limits_{d}F(d).
\end{equation}
Now, further simplifying the expression of $N_5(p),$ by using \eqref{2.18} and Lemma \ref{lem-2.1} we deduce that
\begin{align}\label{N(p)}
 N_5(p)=&\sum\limits_{c,d}\{F(d)+(p-2) \}\{F(c)+(p-2) \}\times \{\phi((c+d)^2-1)+1 \}\notag \\
=&\sum\limits_{c,d}\bigg(\phi((c+d)^2-1)F(c)F(d)+(p-2)^2\phi((c+d)^2-1) +F(c)F(d) \notag\\
&+(p-2)(F(c)+F(d))\phi((c+d)^2-1) + (p-2)^2+(p-2)(F(c)+F(d)) \bigg).
\end{align}
Applying Lemma \ref{lem-2.1} in the following expression, we have  
\begin{align}\label{id-1}
\sum\limits_{c,d}\phi((c+d)^2-1)=\sum\limits_{c}\sum\limits_{d}\phi((c+d)^2-1)=\sum\limits_{c}(-1)=-p.
\end{align}
Moreover, applying Lemma \ref{lem-2.1} and \eqref{2.18} we obtain,
\begin{align}\label{id-2}
\sum\limits_{c,d}\phi((c+d)^2-1)F(c)=\sum\limits_{c} F(c)\sum\limits_{d}\phi((c+d)^2-1)=\sum\limits_{c}-F(c)=-1.
\end{align}
Again, by using \eqref{2.18} we have
\begin{align}\label{id-3}
\sum\limits_{c,d}F(c)F(d)=\sum\limits_{c}F(c)\sum\limits_{d}F(d)=1,
\end{align}
and
\begin{align}\label{id-4}
\sum\limits_{c,d}(F(c)+F(d))=2p.
\end{align}
Finally, substituting \eqref{id-1}, \eqref{id-2}, \eqref{id-3}, and \eqref{id-4} into 
\eqref{N(p)} we obtain 
\begin{align*}
N_5(p)=S+p^4-5p^3+10p^2-10p+5.
\end{align*}

\end{proof}
To this end, for $\lambda\in\F_p,$ we define ${_2F_1}(\lambda):={_2F_1}\left( \begin{matrix}
\phi,&\phi \\
&\varepsilon
\end{matrix}
\middle| \lambda\right).$ We use this notation throughout the paper for brevity.
\begin{lemma}\label{lem-2}
For $t\in\F_p,$ suppose that $F(t)$ is same as defined in Proposition \ref{proposition-1}. Then we have 
    $$ F(t)=p\cdot  {_2 F_1}\left(\frac{t^2}{4}\right)-1 +(p-1)\delta(t),$$
           where $\delta(t)$ is defined by
           $$\delta(t):=\begin{cases}
      1, \ & \text{if} \ t=0;\\
      0, \ & \text{otherwise.}
  \end{cases}$$
  \end{lemma}

\begin{proof} 
For $t=0,$ the identity holds trivially. For $t=2,$ we write
\begin{align*}
           F(2)&=\displaystyle \sum\limits_{x} \phi(x^2-1)\phi((x+2)^2-1)
           =\sum\limits_{x\neq-1} \phi(x-1)\phi(x+3) \\
           &=\sum\limits_{x\neq-1} \phi((x+1)^2-4) 
            =-1-\phi(-1).
            \end{align*}
            Note that the last equality is obtained by using Lemma \ref{lem-2.1}.
Proceeding similarly, we have $F(-2)=-1-\phi(-1).$ Then using \eqref{2F1} and  \cite[Eq. (2.12)]{Greene}, we verify the identity for $t=\pm2$.
Now, suppose that $t\neq0,\pm2.$ Then making the transformation $x\mapsto\frac{1+y}{1-y}$ in the expression of $F(t),$ we have 
    \begin{align*}
    F(t)&=\displaystyle \sum\limits_{y\neq1}\phi\left(\frac{4y}{(1-y)^2}\right)\phi\left(\frac{2+t-ty}{1-y}\right)\phi\left(\frac{y(2-t)+t}{1-y}\right)\\
    &=\sum\limits_{y}\phi(y)\phi\left(\frac{2+t}{t}-y\right)\phi\left(\frac{(2-t)y}{t}+1\right)-1 =p \cdot{_2F_1}\left(\frac{t^2}{4}\right) -1.
    \end{align*}
    Note that the last equality is obtained by transforming $y$ by $\frac{(2+t)y}{t}$ and then using Definition \ref{hypergeometric-2} and \cite[Theorem 4.4(i), Theorem 4.2(ii)]{Greene}.
   
\end{proof}

\begin{proposition}\label{lem-0}
Let $p>3$ be a prime and $S=\sum\limits_{c,d}\phi((c+d)^2-1)F(c)F(d),$ where $F(t)$ is same as defined in Proposition \ref{proposition-1}. Then we have
    $$S=p^2 S_0 +2p(p-1) S_1+(p+2)+(p-1)^2\phi(-1),$$ where $S_0:=\sum\limits_{t}
\phi(t^2-1)
\sum\limits_{s}
{_2F_1}\left(
\frac{s^2}{4}
\right){_2F_1}\left(\frac{(t-s)^2}{4}
\right),$ and $S_1:=\sum\limits_{t}
\phi(t^2-1) {_2F_1}\left(\frac{t^2}{4}
\right).$ 
 
\end{proposition}

\begin{proof} We have
$$
S= \sum\limits_{c,d}\phi((c+d)^2-1)F(c)F(d)
=
\sum\limits_{t}
\phi(t^2-1)
\sum\limits_{s} F(s)F(t-s).
$$
Then using Lemma \ref{lem-2}, for any fixed $t\in\F_p,$ we have

\begin{align}\label{eqn-1000}
\sum\limits_{s} F(s)F(t-s)
=\sum\limits_{s}
\left(
p\cdot{_2F_1}\left(\frac{s^2}{4}\right)
-1+(p-1)\delta(s)\right)
\left(p\cdot{_2F_1}\left(\frac{(t-s)^2}{4}\right)
-1+(p-1)\delta(t-s)\right).
\end{align}
It is routine to verify that
\begin{align*}
    \sum\limits_{s} {_2F_1}\left(\frac{s^2}{4}\right)&=\frac{2}{p}.
 \end{align*}
Now, substituting the above identity into \eqref{eqn-1000},
we deduce that
\begin{align*}
\sum\limits_{s} F(s)F(t-s)
=&p^2
\sum\limits_{s} {_2F_1}\left(
\frac{s^2}{4}
\right){_2F_1}\left(\frac{(t-s)^2}{4}
\right)
+
2p(p-1)\ {_2F_1}\left(
\frac{t^2}{4}
\right)-(p+2)
+(p-1)^2\delta(t).
\end{align*}
Finally, substituting the above identity in the expression of $S$ and simplifying further we arrive at the required result.

\end{proof}

\begin{lemma}\label{lem-4}
Let $p\geq5$ be a prime and $S_1:=\sum\limits_{t}
\phi(t^2-1) {_2F_1}\left(\frac{t^2}{4}
\right).$ Then we have 
   $$ S_1=\frac{1+\phi(-1)}{p}+\phi(-3)+p\cdot  {_3F_2}\left( \begin{matrix}
\phi,&\phi,&\phi \\
&\varepsilon,&\varepsilon
\end{matrix}
\middle| \frac{1}{4}\right).$$
\end{lemma}

\begin{proof}

Recall that
    \begin{align}  S_1&=\sum\limits_{x}\phi(x^2-1)\cdot {_2F_1}\left(\frac{x^2}{4}\right)=\sum\limits_{x}\left(1+\phi(x)\right)\phi(x-1){_2F_1}\left( \frac{x}{4}\right)\notag\\
    &=\sum\limits_{x}\phi(x-1){_2F_1}\left( \frac{x}{4}\right)+\sum\limits_{x}\phi(x)\phi(x-1){_2F_1}\left( \frac{x}{4}\right)\label{eqn-2001}.
\end{align}
Now, simplifying the first summation of that arise in $S_1$ and using \cite[Eq. (2.6)]{Greene}, we obtain that
\begin{align}\label{eqn-2000}
\sum\limits_{x}\phi(x-1)\ {_2F_1}\left(\frac{x}{4}\right)&=\frac{p}{p-1}\sum\limits_{x}\phi(x-1)\sum\limits_{\chi}\binom{\phi\chi}{\chi}^2\chi\left(\frac{x}{4}\right)\notag \\
&=\frac{p}{p-1}\sum\limits_{\chi}\binom{\phi\chi}{\chi}^2\overline{\chi}(4)\sum\limits_{x}\phi(x-1)\chi(x)\notag\\
&=\frac{p^2}{p-1}\sum\limits_{\chi}\binom{\phi\chi}{\chi}^2\ \binom{\chi}{\phi\chi}\overline{\chi}(4).
\end{align}
Using \cite[Eq (2.12), Eq (2.13), Eq (2.10)]{Greene} and Lemma \ref{lem-3} in \eqref{eqn-2000}, we simplify the above expression further and derive that
\begin{align}
\sum\limits_{x}\phi(x-1)\ {_2F_1}\left(\frac{x}{4}\right)
=\frac{1+\phi(-1)}{p}+\phi(-3).\label{eqn-2002}
\end{align}
Applying \cite[Theorem 3.13]{Greene} in the following expression, we have
\begin{align}
  \sum\limits_{x}\phi(x)\phi(x-1)\ {_2F_1}\left( \frac{x}{4}\right)&=p\cdot{_3F_2}\left( \begin{matrix}
\phi,&\phi,&\phi \\
&\varepsilon,&\varepsilon
\end{matrix}
\middle| \frac{1}{4}\right)\label{eqn-2003}.
\end{align}
Finally, substituting \eqref{eqn-2002} and \eqref{eqn-2003} into \eqref{eqn-2001}, we deduce the result.
\end{proof}

\begin{lemma}\label{lem-5}
Let $p\geq5$ be a prime and $S_0=\sum\limits_{t}
\phi(t^2-1)
\sum\limits_{s}
{_2F_1}\left(
\frac{s^2}{4}
\right){_2F_1}\left(\frac{(t-s)^2}{4}
\right).$ Then 
   \begin{align*}
& S_0=2c_p+\frac{2}{p^2(p-1)}\left(1+\phi(-1)\right)+\frac{2\phi(-1)}{p(p-1)}\left(\phi(-1)+1+(p-1)\phi(3)\right)\\
& \hspace{2cm}+\frac{p^2}{p-1}\ \sum\limits_{\psi\neq\varepsilon,\phi}\overline{\psi}(-4)\ {_4F_3}\left(\begin{matrix}
\overline{\psi},&\phi ,&\phi,&\phi\\
& \phi\overline{\psi},& \phi\overline{\psi},& \phi\overline{\psi}\\
\end{matrix}
\middle|1
\right)\\
& \hspace{2cm}+\frac{p^3\phi(-1)}{p-1}\ \sum\limits_{\psi\neq\varepsilon,\phi}\overline{\psi}(-4)\binom{\psi}{\phi}^2{_4F_3}\left(\begin{matrix}
\overline{\psi},&\phi ,&\phi,&\phi\\
& \phi\overline{\psi},& \phi\overline{\psi},& \phi\overline{\psi}\\
\end{matrix}
\middle|1
\right)+T',
\end{align*}
 where $$  T'=\begin{cases}
-\frac{p^2-3p-2-a_h(p)}{p^2} ,&\text{if} \ p\equiv 3 \pmod{4};\\
\frac{p^3-2p^2-5p-10-(p+1)a_h(p)}{p^2(p-1)} \ , & \text{if} \ p\equiv 1 \pmod{4},
\end{cases}$$ 
and where $a_h(p),$ and $c_p$ are same as defined in \eqref{b_p}, and \eqref{c_p}, respectively.
\end{lemma}

\begin{proof} \allowdisplaybreaks
We have, 
\begin{align*} S_0&=\sum\limits_{t}
\phi(t^2-1)
\sum\limits_{s} {_2F_1}\left(  \frac{s^2}{4}\right){_2F_1}\left(\frac{(t-s)^2}{4}\right)\\
&=\frac{p^2}{(p-1)^2}\sum\limits_{t\neq0}
\phi(t^2-1)\sum\limits_{\chi_1,\chi_2,s}\binom{\phi \chi_1}{\chi_1}^2\binom{\phi \chi_2}{\chi_2}^2\chi_1\left(\frac{s^2}{4}\right)\chi_2\left(\frac{(t-s)^2}{4}\right)+T_0,
\end{align*}
where $T_0:=\frac{p^2\phi(-1)}{(p-1)^2}\sum\limits_{\chi_1,\chi_2,s}\binom{\phi \chi_1}{\chi_1}^2\binom{\phi \chi_2}{\chi_2}^2(\chi_1\chi_2)^2\left(\frac{s}{2}\right).$\\
Now, making the transformation $s\rightarrow ts$ and subsequently $t\rightarrow\frac{1+t}{1-t}$ in the main term of $S_0,$ we deduce that 
\begin{align*}
S_0&=\frac{p^3}{(p-1)^2}\sum\limits_{\chi_1,\chi_2}\binom{\phi \chi_1}{\chi_1}^2\binom{\phi \chi_2}{\chi_2}^2J(\chi_1^2,\chi_2^2) (\overline{\chi_1\chi_2})(4)\ {_2F_1}\left(\begin{matrix}
\left(\overline{\chi_1\chi_2}\right)^2, &\phi \\
&\phi\left(\overline{\chi_1\chi_2}\right)^2\\    
\end{matrix}
\middle|-1\right)+T_0.
\end{align*}
Again, replacing $\chi_1\chi_2$ by $\psi$ and $\chi_1$ by $\chi,$ in the main term of $S_0,$ we write
\begin{align}
S_0&=\frac{p^3}{(p-1)^2} \displaystyle\sum\limits_{\chi,\psi}\binom{\phi \chi}{\chi}^2\binom{\phi \psi\overline{\chi}}{\psi\overline{\chi}}^2J(\chi^2,(\psi\overline{\chi})^2) \overline{\psi}(4)\ {_2F_1}\left(\begin{matrix}
\overline{\psi}^2, & \phi \\
&\phi\overline{\psi}^2 \\    
\end{matrix}
\middle|-1
\right)+T_0\notag\\
&=T'+ \frac{p^3}{(p-1)^2} \displaystyle\sum\limits_{\chi\in \widehat{{\F}_p^{\times}},\psi\neq\{ \varepsilon,\phi \}}\binom{\phi \chi}{\chi}^2\binom{\phi \psi\overline{\chi}}{\psi\overline{\chi}}^2J(\chi^2,(\psi\overline{\chi})^2) \overline{\psi}(4)\ {_2F_1}\left(\begin{matrix}
\overline{\psi}^2, & \phi \\
&\phi\overline{\psi}^2 \\    
\end{matrix}
\middle|-1
\right),\label{2004}
\end{align}
where $ T'=\frac{p^3}{(p-1)^2}\cdot {_2F_1}\left(\begin{matrix}
\varepsilon, & \phi \\
&\phi \\    
\end{matrix}
\middle|-1
\right) \left\{\sum\limits_{\chi}\binom{\phi \chi}{\chi}^2\binom{\phi\overline{\chi}}{\overline{\chi}}^2J(\chi^2,\overline{\chi}^2)+\sum\limits_{\chi}\binom{\phi \chi}{\chi}^2\binom{\overline{\chi}}{\phi\overline{\chi}}^2J(\chi^2,\overline{\chi}^2)\right\}+T_0.$\\
Using \cite[Eq (2.7)]{Greene}, and \eqref{Jacobi-2}, for $\psi\neq\varepsilon,\phi,$ we have 
\begin{align*}
 \binom{\phi \chi}{\chi}^2\binom{\phi \psi\overline{\chi}}{\psi\overline{\chi}}^2J(\chi^2,(\psi\overline{\chi})^2)&= \binom{\phi}{\chi}^2\binom{\phi}{\psi\overline{\chi}}^2 \frac{g(\chi^2)g((\psi\overline{\chi})^2)}{g(\psi^2)}.
 \end{align*}
 Applying \eqref{Hasse}, \eqref{Gauss}, \cite[Eq. (2.8)]{Greene}, and \eqref{Jacobi-2}, in above expression, we derive  
 \begin{align*}
 \binom{\phi \chi}{\chi}^2\binom{\phi \psi\overline{\chi}}{\psi\overline{\chi}}^2J(\chi^2,(\psi\overline{\chi})^2)
 &=\frac{p^2}{\phi(-1)}\binom{\psi}{\phi}\binom{\overline{\psi}\chi}{\chi}^2\binom{\phi\chi}{\psi}^3\binom{\chi}{\overline{\psi}\chi}.
 \end{align*}
Substituting the above expression in \eqref{2004}, and then using \cite[Eq (2.6)]{Greene}, we write
 \begin{align*}
    S_0&=T'+ \frac{p^5\phi(-1)}{(p-1)^2} \sum\limits_{\psi\neq\{ \varepsilon,\phi \}} \overline{\psi}(4)\ \binom{\psi}{\phi}\left\{ \frac{p}{p-1}\sum\limits_{\chi_1}\binom{\overline{\psi}^2\chi_1}{\chi_1}\binom{\phi\chi_1}{\phi\overline{\psi}^2\chi_1}\chi_1(-1)\right\} \\
& \times\left\{  \binom{\overline{\psi}}{\varepsilon}^2\binom{\phi}{\psi}^3\binom{\varepsilon}{\overline{\psi}}+\binom{\varepsilon}{\psi}^2\binom{\phi\psi}{\psi}^3\binom{\psi}{\varepsilon}+\sum\limits_{\chi\neq\{\varepsilon,\psi\}}\binom{\overline{\psi}\chi}{\chi}^2\binom{\phi\chi}{\phi\overline{\psi}\chi}^3\binom{\chi}{\overline{\psi}\chi}\right\}.
\end{align*}
Now, in the above expression of $S_0,$ using \cite[Eq (2.6), Eq (2.7), Eq (2.8), Eq (2.12), Eq (2.13)]{Greene} we obtain that
\begin{align*}
S_0&=T'+ \frac{p^5\phi(-1)}{(p-1)^2} \sum\limits_{\psi\neq\{ \varepsilon,\phi \}} \overline{\psi}(4)\ \binom{\psi}{\phi}\left\{ \frac{p}{p-1}\sum\limits_{\chi_1} \binom{\psi^2}{\chi_1}\chi_1(-1)\binom{\phi\chi_1}{\psi^2}\chi_1(-1) \right\} \\
&\times\left\{  -\frac{2\phi(-1)}{p^3}\binom{\overline{\psi}}{\phi}^3+\sum\limits_{\chi\neq\{\varepsilon,\psi\}}\binom{\overline{\psi}\chi}{\chi}^2\binom{\phi\chi}{\phi\overline{\psi}\chi}^3\binom{\chi}{\overline{\psi}\chi}\right\}.
\end{align*}
Simplifying the above expression by applying \cite[Eq (2.15), Eq (2.8)]{Greene}, and Lemma \ref{lem-3}, we have
\begin{align*}
S_0&=T'+ \frac{p^5\phi(-1)}{(p-1)^2} \sum\limits_{\psi\neq\{ \varepsilon,\phi \}} \overline{\psi}(4)\ \binom{\psi}{\phi}\left\{ \frac{p}{p-1}\sum\limits_{\chi_1}\binom{\phi\chi_1}{\chi_1}\binom{\overline{\psi}^2\chi_1}{\phi}\phi\overline{\chi_1}(-1)\right\}  \\
& \times \left\{ -\frac{2\phi(-1)}{p^3}\binom{\overline{\psi}}{\phi}^3 + \frac{\psi(-1)}{p}\sum\limits_{\chi}\binom{\overline{\psi}\chi}{\chi}\binom{\phi\chi}{\phi\overline{\psi}\chi}^3-\frac{\psi(-1)}{p}\left(\binom{\overline{\psi}}{\varepsilon}\binom{\phi}{\phi\overline{\psi}}^3+\binom{\varepsilon}{\psi}\binom{\phi\psi}{\phi}^3 \right) \right\} \\
&=T'+ \frac{p^5\phi(-1)}{(p-1)^2} \sum\limits_{\psi\neq\{ \varepsilon,\phi \}} \overline{\psi}(4)\ \binom{\psi}{\phi}\left\{ \frac{p}{p-1}\sum\limits_{\chi_1}\binom{\phi\chi_1}{\chi_1}\binom{\overline{\psi}^2\chi_1}{\phi\overline{\psi}^2\chi_1}\phi\chi_1(-1) \right\} \\
& \times \left\{ \frac{2(p-1)\phi(-1)}{p^3}\binom{\overline{\psi}}{\phi}^3+\frac{(p-1)\psi(-1)}{p^2} {_4F_3}\left(\begin{matrix}
\overline{\psi},&\phi ,&\phi,&\phi\\
& \phi\overline{\psi},& \phi\overline{\psi},& \phi\overline{\psi}\\
\end{matrix}
\middle|1
\right) \right\}.
\end{align*}
Note that the last equality is obtained by using \cite[Eq (2.6), Eq (2.7), Eq (2.8), Eq (2.12), Eq (2.13)]{Greene}.
Now, splitting all the terms in the above expression, we derive

\begin{align}\label{2005}
S_0&=T'+ \frac{p^2\phi(-1)}{(p-1)} \sum\limits_{\psi\neq\{ \varepsilon,\phi \}}\bigg\{\frac{2\overline{\psi}(4)}{p}\binom{\overline{\psi}}{\phi}^3+2\overline{\psi}(4)\phi(-1)\binom{\overline{\psi}}{\phi}^3\binom{\psi}{\phi}^2 \notag\\
&+\overline{\psi}(4)\phi\psi(-1)\ {_4F_3}\left(\begin{matrix}
\overline{\psi},&\phi ,&\phi,&\phi\\
& \phi\overline{\psi},& \phi\overline{\psi},& \phi\overline{\psi}\\
\end{matrix}
\middle|1
\right)
+p\ \overline{\psi}(4)\psi(-1)\ {_4F_3}\left(\begin{matrix}
\overline{\psi},&\phi ,&\phi,&\phi\\
& \phi\overline{\psi},& \phi\overline{\psi},& \phi\overline{\psi}\\
\end{matrix}
\middle|1
\right)
\binom{\psi}{\phi}^2\bigg\}.
\end{align}
Now using \cite[Eq (2.7), Eq (2.8)]{Greene}, we obtain 
\begin{align}\label{2006}
    \sum\limits_{\psi\neq\varepsilon,\phi}\overline{\psi}(4)\frac{2}{p}\binom{\overline{\psi}}{\phi}^3&=\frac{2}{p^3}\sum\limits_{\psi\neq\varepsilon,\phi}\overline{\psi}(4)\binom{\phi}{\psi}^3\phi\psi(-1)
    =\frac{2}{p}\sum\limits_{\psi\neq\varepsilon,\phi}\overline{\psi}(4)\binom{\phi\psi}{\psi}^3\phi(-1)\notag\\
    &=\frac{2\phi(-1)}{p}\cdot \frac{p-1}{p} {_3F_2}\left(\begin{matrix}
\phi ,&\phi,&\phi\\
& \varepsilon,& \varepsilon\\
\end{matrix}
\middle|\frac{1}{4}
\right)-\frac{2\phi(-1)}{p}\left[\binom{\phi}{\varepsilon}^3+\binom{\varepsilon}{\phi}^3\right]\notag\\
&=\frac{2\phi(-1)(p-1)c_p}{p^2}+\frac{2\phi(-1)}{p^4}\left(1+\phi(-1)\right).
    \end{align}
    Note that the last equality is obtained by applying \cite[Eq (2.12), Eq (2.13)]{Greene} and \eqref{3F2}. Now, applying Lemma \ref{lem-3} and \cite[Eq (2.5)]{Greene}, it is straightforward to verify that
    \begin{equation}\label{2200}\sum\limits_{\psi\neq\varepsilon,\phi}2\overline{\psi}(4)\phi(-1)\binom{\overline{\psi}}{\phi}^3\binom{\psi}{\phi}^2=\frac{2}{p^3}\left(\phi(-1)+1+(p-1)\phi(3)\right).\end{equation}
    Substituting \eqref{2006} and \eqref{2200} into \eqref{2005}, we obtain that
    \begin{align}
& S_0=2\ c_p+\frac{2}{p^2(p-1)}\left(1+\phi(-1)\right)+\frac{2\phi(-1)}{p(p-1)}\left(\phi(-1)+1+(p-1)\phi(3)\right)\\
& \hspace{2cm}+\frac{p^2}{p-1}\ \sum\limits_{\psi\neq\varepsilon,\phi}\overline{\psi}(-4)\ {_4F_3}\left(\begin{matrix}
\overline{\psi},&\phi ,&\phi,&\phi\\
& \phi\overline{\psi},& \phi\overline{\psi},& \phi\overline{\psi}\notag\\
\end{matrix}
\middle|1
\right)\notag\\
& \hspace{2cm}+\frac{p^3\phi(-1)}{p-1}\ \sum\limits_{\psi\neq\varepsilon,\phi}\overline{\psi}(-4)\binom{\psi}{\phi}^2{_4F_3}\left(\begin{matrix}
\overline{\psi},&\phi ,&\phi,&\phi\\
& \phi\overline{\psi},& \phi\overline{\psi},& \phi\overline{\psi}\\
\end{matrix}
\middle|1
\right)+T'.\notag
\end{align}
Now, we evaluate the expression for $T'.$ Following \cite[Eq (4.11), Eq (2.12)]{Greene} we have
\begin{align}\label{2009} 
{_2F_1}\left(\begin{matrix}
\varepsilon, & \phi \\
&\phi \\    
\end{matrix}
\middle|-1
\right)
&=\begin{cases}
0, & p\equiv 3 \pmod{4};\\
 \frac{-2}{p}, & p\equiv1 \pmod{4}.
\end{cases}
\end{align}
Using \eqref{Gauss},\eqref{Jacobi-2} and $g(\varepsilon)=-1$ , for $\chi\neq\varepsilon,\phi,$ we have
$$J(\chi^2,\overline{\chi}^2)
=\begin{cases}
p-2, & \text{if}\ \chi=\varepsilon,\phi; \\
-1, & \text{Otherwise}.
\end{cases}$$
Using this value and \cite[Eq (2.12), Eq (2.13)]{Greene}, we write
\begin{align*}
& \sum\limits_{\chi}\binom{\phi \chi}{\chi}^2\binom{\phi\overline{\chi}}{\overline{\chi}}^2J(\chi^2,\overline{\chi}^2)
= \frac{2(p-2)}{p^4}-\sum\limits_{\chi\neq \{\varepsilon,\phi\}}\binom{\phi \chi}{\chi}^2\binom{\phi\overline{\chi}}{\overline{\chi}}^2.
\end{align*}
Following \cite[Eq (2.8)]{Greene} and Lemma \ref{lem-3} we have
\begin{align}\label{2008}
\sum\limits_{\chi}\binom{\phi \chi}{\chi}^2\binom{\phi\overline{\chi}}{\overline{\chi}}^2J(\chi^2,\overline{\chi}^2)&=\frac{2(p-2)}{p^4}-\frac{p-3}{p^2}=-\frac{(p-1)(p^2-2p-4)}{p^4}.
\end{align}
Using similar calculations, we derive that
\begin{align}\label{2007}
\sum\limits_{\chi}\binom{\phi \chi}{\chi}^2\binom{\overline{\chi}}{\phi\overline{\chi}}^2J(\chi^2,\overline{\chi}^2)
&=\frac{2(p-1)}{p^4} - \frac{p-1}{p} {_4F_3} 
    \left(\begin{matrix}
    \phi,&\phi,&\phi,&\phi\\
    &\varepsilon,&\varepsilon,&\varepsilon  \\   
\end{matrix}
\middle|1
\right) \notag \\ 
&=\frac{p-1}{p}\cdot \frac{a_h(p)+p}{p^3}+\frac{2(p-1)}{p^4},
\end{align}
and 
\begin{align}
    \sum\limits_{\chi_1,\chi_2,s}\binom{\phi \chi_1}{\chi_1}^2\binom{\phi \chi_2}{\chi_2}^2(\chi_1\chi_2)^2\left(\frac{s}{2}\right)
    &=\frac{(p-1)^2(p^2-2p-2)}{p^4}+\frac{p-1}{p} {_4F_3} 
    \left(\begin{matrix}
    \phi,&\phi,&\phi,&\phi\\
    &\varepsilon,&\varepsilon,&\varepsilon \\   
\end{matrix}
\middle|1
\right) \notag \\
&=\frac{(p-1)^2}{p^4}(p^2-3p-2-a_h(p)), \label{130}
\end{align}
where $a_h(p)$ is as defined in \eqref{b_p}.
Note that the last equalities in the above two expressions are due to \cite[Theorem 6]{AO}.
Substituting \eqref{2009}, \eqref{2008} , \eqref{2007} and \eqref{130} in the expression of $T'$ we complete the proof.

\end{proof}

\begin{proposition}\label{prop-1}
For a prime $p>5,$ we have 
$$\sum\limits_{\psi\neq \varepsilon,\phi}\overline{\psi}(-4){_4F_3}\left(\begin{matrix}
\overline{\psi},&\phi ,&\phi,&\phi\\
& \phi\overline{\psi},& \phi\overline{\psi},& \phi\overline{\psi}\\
\end{matrix}
\middle|1
\right)=\frac{p-1}{p^3}(2p\phi(-3)+\frac{p+2}{p-1}+a_g(p)+\frac{\phi(-1)}{p-1}(a_h(p)+p)),$$
where $a_g(p)$ and $a_h(p)$ are same as defined in \eqref{a_p} and \eqref{b_p}, respectively.
    
\end{proposition}

\begin{proof}
    Recall the classical Kloosterman sum $K(a,b)$ over a finite field $\F_p$ is given by $$K(a,b)=\sum\limits_{x\in\F_p^{\times}}\zeta_p^{ax+bx^{-1}}.$$ For our purpose, if we choose $b=1,$ then we have $K(a,1)=K(a).$ Now, let $$M_5(p):=\#\{(x_1,x_2,x_3,x_4,x_5)\in(\F_p^{\times})^5:\sum\limits_{x=1}^{5}x_i=0\  \text{and}\ \sum\limits_{x=1}^{5}x_i^{-1}=0 \}.$$ Then using \cite[Section 4.4]{iwaniec}, we have 
    \begin{equation}\label{S_5}(p-1)S(5)=M_5(p) p^2-(p-1)^5+2(p-1).
    \end{equation}
     
    Now, substituting $x_1=x_1x_5,x_2=x_2x_5,x_3=x_3x_5,x_4=x_4x_5,x_5=x_5$ in $M_5(p),$ we have $$M_5(p)=(p-1)\cdot \#\{(x_1,x_2,x_3,x_4)\in(\F_p^{\times})^4:\sum\limits_{x=1}^{4}x_i=-1\  \text{and}\ \sum\limits_{x=1}^{4}x_i^{-1}=-1\}.$$
 
 For $a,b\in\F_p$ and $k\ge1$, define $N_{k,p}(a,b):=\# \{(x_1,x_2,\cdots,x_k)\in(\F_p^{\times})^k:\sum\limits_{x=1}^{k}x_i=a\  \text{and}\ \sum\limits_{x=1}^{k}x_i^{-1}=b\}.$ Using this notation, we have that 
 
 \begin{equation}\label{M_5}
    M_5(p)=(p-1)\cdot N_{4,p}(-1,-1).
    \end{equation}
 
 It is trivial that
    $$N_{k,p}(-a,-b)=N_{k,p}(a,b).$$
    Therefore, we have 
\begin{equation}\label{N_4} N_{4,p}(-1,-1)=\sum\limits_{u,v\in\F_p}N_{2,p}(u,v)N_{2,p}(1+u,1+v).\end{equation}
   Observe that in the expression of $N_{2,p}(u,v)$, we have $u=0$ if and only if $v=0$. Therefore, using this fact in \eqref{N_4} we obtain that
     \begin{equation}\label{N_4(-1,-1)}  N_{4,p}(-1,-1)=\sum\limits_{u,v\neq{0,-1}}N_{2,p}(u,v)N_{2,p}(1+u,1+v)+2N_{2,p}(0,0)N_{2,p}(1,1).\end{equation}
     For $u,v\neq0,$ we have the following: \begin{equation}\label{d} N_{2,p}(u,v)=1+\phi\left(1-\frac{4}{uv}\right), N_{2,p}(0,0)=p-1, \ \text{and } \ N_{2,p}(1,1)=1+\phi(-3).\end{equation}
     Now, \eqref{d} yields, \begin{align}
        & \sum\limits_{u,v\neq{0,-1}}N_{2,p}(u,v)N_{2,p}(1+u,1+v)\notag\\ &=\sum\limits_{u,v\neq{0,-1}}\left\{1+\phi\left(1-\frac{4}{uv}\right)\right\}\left\{1+\phi\left(1-\frac{4}{(u+1)(v+1)}\right)\right\} \notag\\
        &=(p-2)^2+ \sum\limits_{u,v\neq{0,-1}}\phi\left(1-\frac{4}{uv}\right)+ \sum\limits_{u,v\neq{0,-1}}\phi\left(1-\frac{4}{(u+1)(v+1)}\right)  +D\label{e},
     \end{align}
     where $D=\sum\limits_{u,v\neq{0,-1}}\phi\left(1-\frac{4}{uv}\right)\phi\left(1-\frac{4}{(u+1)(v+1)}\right).$\\
     It is easy to verify that $\sum\limits_{u,v\neq{0,-1}}\phi\left(1-\frac{4}{uv}\right)=\sum\limits_{u,v\neq{0,-1}}\phi\left(1-\frac{4}{(u+1)(v+1)}\right)=\phi(-3)-(p-3).$\\
     Now, substituting $u=\frac{x}{1-x}$ and $v=\frac{y}{1-y}$ in $D$ we obtain
     $$D=\sum\limits_{x,y\neq0,1}\phi(1-4(1-x)(1-y))\phi(1-4(1-x^{-1})(1-y^{-1})).$$
     Let $A:=\sum\limits_{\psi\in\widehat{\F_p^{\times}}}\overline{\psi}(-4){_4F_3}\left(\begin{matrix}
\overline{\psi},&\phi ,&\phi,&\phi\\
& \phi\overline{\psi},& \phi\overline{\psi},& \phi\overline{\psi}\\
\end{matrix}
\middle|1
\right).$
     Then, we have \begin{equation}\label{z}
        A=\frac{1}{p^3}\sum\limits_{\substack{x,y,z\neq0,1 \\ xyz\neq1}} \phi(xyz)\sum\limits_{\psi}\psi \left(\frac{4(1-x)(1-y)(1-z)}{1-xyz} \right) .   
\end{equation}
Using \eqref{ortho-2}, we have that the inner sum in \eqref{z} is nonzero only if $$
  \frac{4(1-x)(1-y)(1-z)}{1-xyz}=1. 
  $$ This implies $z=\frac{1-4(1-x)(1-y)}{xy-4(1-x)(1-y)}.$
  Substituting the value of $z$ into \eqref{z}, we obtain \begin{align} A
&=\frac{p-1}{p^3}\sum\limits_{\ x\neq0,1} \left\{ \sum\limits_{y\neq 0,1}\phi\left(xy\frac{1-4(1-x)(1-y)}{xy-4(1-x)(1-y)}\right)-1\right\}\ \notag\\
&= \frac{p-1}{p^3}\sum\limits_{\ x,y\neq0,1} \phi\left(xy\frac{1-4(1-x)(1-y)}{xy-4(1-x)(1-y)}\right)-\frac{(p-1)(p-2)}{p^3}\notag\\
&=\frac{p-1}{p^3}(D-(p-2)).\label{b}
\end{align}
Now, substituting \eqref{d}, and \eqref{e} into \eqref{N_4(-1,-1)}, we obtain 
\begin{align}N_{4,p}(-1,-1)&=(p-2)^2+2(\phi(-3-(p-3)))+D+2(p-1)(1+\phi(-3))\notag\\
&= D+2p\phi(-3)+p^2-4p+8.\label{a}\\ 
\notag\end{align}
Using \eqref{M_5},\eqref{a}, and \eqref{b} into \eqref{S_5} we obtain 
\begin{equation}\label{c}
S(5)=\frac{p^5}{p-1}A+p^3+2p^3\phi(-3)+4p+1.\end{equation}
Observe that for a prime $p>3$, we have   $\leg{p}{3}=\phi(-3).$
Now substituting \eqref{s(5)} into \eqref{c} we obtain 
\begin{equation}\label{f}A=\frac{p-1}{p^3}(2p\phi(-3)-p+a_g(p)+5).\end{equation}
Using \eqref{f}, and \cite[Eq (2.12), Eq (2.14) ]{Greene}, we conclude the proof. 
\end{proof}

The following result gives an explicit expression for $N_k(p)$ for $k\geq1$, providing a connection between the point counting problem of $X_k$ and $k$-th moments of Kloosterman sums and its quadratic twist.

\begin{proposition}\label{theorem-2}
    For prime $p\geq3,$ and $k\geq1,$ we have  $$ N_k(p)=\frac{(p-1)^k}{p}+\frac{1}{p}S(k)+\frac{1}{p}S(k,\phi).$$ 

\end{proposition}

\begin{proof}
 Recall that
 $ N_k(p)=\#\left\{(x_1,x_2,\cdots,x_k)\in({\F}_p^{\times})^k:
\sum\limits_{i=1}^{k}(x_i+x_i^{-1})=0\right  \} .$
Now, using the fact that for $t\in\F_p,$\begin{align*}\sum\limits_{x\in \F_p}\zeta_p^{xt} =
\begin{cases}
p, &\text{if}\  t=0 ;\\
0, &\text{if}\  t\neq 0,\\
\end{cases}
\end{align*}
we have
\begin{align*}N_k(p)
 &=\frac{1}{p}\sum\limits_{t}\sum\limits_{x_1,x_2,\cdots,x_k\in{\F}_p^{\times}}\zeta_p^{t\left(\sum\limits_{i=1}^{k}(x_i+x_i^{-1})\right)} \\
&=\frac{(p-1)^k}{p}+\frac{1}{p}\sum\limits_{t,x_1,x_2,\cdots,x_k\in{\F}_p^{\times}}\zeta_p^{t\left(\sum\limits_{i=1}^{k}(x_i+x_i^{-1})\right) }\\
&=\frac{(p-1)^k}{p}+\frac{1}{p}\sum\limits_{t,x_1,x_2,\cdots,x_k\in{\F}_p^{\times}}\zeta_p^{\sum\limits_{i=1}^{k}(x_i+t^2x_i^{-1})}\\
&=\frac{(p-1)^k}{p}+\frac{1}{p}\sum\limits_{t,x_1,x_2,\cdots,x_k\in{\F}_p^{\times}}\prod_{i=1}^{k}\zeta_p^{(x_i+t^2x_i^{-1})} \\
&= \frac{(p-1)^k}{p}+\frac{1}{p}\sum\limits_{t\neq0}\prod_{i=1}^{k}\left(\sum\limits_{x_i\neq0}\zeta_p^{(x_i+t^2x_i^{-1})}\right) \\
&= \frac{(p-1)^k}{p}+\frac{1}{p}\sum\limits_{t\neq0}K(t^2)^k \\
&=\frac{(p-1)^k}{p}+\frac{1}{p}\sum\limits_{t\neq0}K(t)^k+\frac{1}{p}\sum\limits_{t\neq0}\phi(t)K(t)^k. 
\end{align*} 
Finally, using \eqref{132}, and \eqref{Twisted-Moments}, we complete the proof.

\end{proof}

\section{Proof of the main theorems}
In this section, we prove our main results. We begin with the proof of Theorem \ref{theorem-1}.

\begin{proof}[Proof of Theorem \ref{theorem-1}]
Recall that $S:=\sum\limits_{c,d}\phi((c+d)^2-1)F(c)F(d),$ where $F(t)$ is same as defined in Proposition \ref{proposition-1}.
Now, using Proposition \ref{lem-0}, Lemma \ref{lem-4}, Lemma \ref{lem-5}, and Proposition \ref{prop-1}, we obtain \begin{align}S&=2p^3 \cdot {_3F_2}\left(\begin{matrix}
        \phi,&\phi, & \phi\\
        &\varepsilon, & \varepsilon
    \end{matrix}
    \middle|\frac{1}{4}
    \right)+pa_g(p)+4p^2\phi(-3)+\frac{(p+2)(2p-1)}{p-1}\notag\\&+\frac{p^5\phi(-1)}{p-1}\sum\limits_{\psi\neq\varepsilon,\phi}\overline{\psi}(-4){_4F_3}\left(\begin{matrix}
\overline{\psi},&\phi ,&\phi,&\phi\\
& \phi\overline{\psi},& \phi\overline{\psi},& \phi\overline{\psi}\\
\end{matrix}
\middle|1
\right)\binom{\psi}{\phi}^2+T \label{k},
\end{align}
where $T$ is same as defined in Theorem \ref{theorem-1}.
From \cite[Eq (17)]{Kalita}, we recall that for any $\chi\in \widehat{{\F}_p^{\times}}\setminus\{\varepsilon\},$
\begin{align}
    S(4,\chi)&=pS(2,\chi)+\phi(-1)p^2g(\chi)^2\overline{\chi}(4)\left\{p\cdot{_4F_3}\left(\begin{matrix}
\chi,&\chi ,&\chi,&\phi\\
& \phi\chi,& \varepsilon,& \varepsilon\\
\end{matrix}
\middle|1
\right)+\frac{1}{p}\binom{\chi}{\phi}\right\}\notag\\&\hspace{2cm}+\phi(-1)p(p-1)g(\chi)^2\overline{\chi}(4)\binom{\chi}{\phi}.\label{g}
\end{align}
Now, substituting \eqref{eqn-Liu} into \eqref{g}, we have
\begin{align}
\phi(-1)\cdot p^3\cdot \overline{\chi}(4)\cdot{_4F_3}\left(\begin{matrix}
\chi,&\chi ,&\chi,&\phi\\
& \phi\chi,& \varepsilon,& \varepsilon\\
\end{matrix}
\middle|1
\right) =\frac{S(4,\chi)-2pS(2,\chi)}{g(\chi)^2}\label{i}.
\end{align}
For $\chi\neq\varepsilon,\phi$, using \cite[Eq (2.7), Eq (2.6), Eq (2.15)]{Greene} we obatain
$${_4F_3}\left(\begin{matrix}
\overline{\chi},&\phi ,&\phi,&\phi\\
& \phi\overline{\chi},& \phi\overline{\chi},& \phi\overline{\chi}\\
\end{matrix}
\middle|1
\right) =\frac{p}{p-1}\sum\limits_{\psi}\binom{\phi\psi}{\psi}\binom{\phi}{\chi\overline{\psi}}\binom{\phi\psi}{\chi}^2\psi(-1).$$ 
Substituting $\psi$ by $\phi\overline{\psi}$ in above expression, we obtain
\begin{align}{_4F_3}\left(\begin{matrix}
\overline{\chi},&\phi ,&\phi,&\phi\\
& \phi\overline{\chi},& \phi\overline{\chi},& \phi\overline{\chi}\\
\end{matrix}
\middle|1
\right)&= \frac{p}{p-1}\sum\limits_{\psi}\binom{\overline{\psi}}{\phi\overline{\psi}}\binom{\phi}{\phi\chi\psi}\binom{\overline{\psi}}{\chi}^2\phi\overline{\psi}(-1).\notag\end{align}
Using \cite[Eq (2.6),Eq (2.7),Eq (2.8)]{Greene} in above expression we obtain,
\begin{align}{_4F_3}\left(\begin{matrix}
\overline{\chi},&\phi ,&\phi,&\phi\\
& \phi\overline{\chi},& \phi\overline{\chi},& \phi\overline{\chi}\\
\end{matrix}
\middle|1
\right)&=\frac{p}{p-1}\sum\limits_{\psi}\binom{\phi\psi}{\psi}\binom{\chi\psi}{\phi\chi\psi}\binom{\overline{\psi}}{\overline{\chi}\overline{\psi}}^2\phi\overline{\psi}(-1)\phi(-1)\phi\chi\psi(-1).\notag
\end{align}
Now, using \cite[Eq (2.8)]{Greene} in above equation, we finally obtain
\begin{align}
{_4F_3}\left(\begin{matrix}
\overline{\chi},&\phi ,&\phi,&\phi\\
& \phi\overline{\chi},& \phi\overline{\chi},& \phi\overline{\chi}\\
\end{matrix}
\middle|1
\right)&=\phi\chi(-1){_4F_3}\left(\begin{matrix}
\chi,&\chi ,&\chi,&\phi\\
& \phi\chi,& \varepsilon,& \varepsilon\\
\end{matrix}
\middle|1
\right).\label{h}
\end{align}
Therefore, for $\psi\neq\varepsilon,\phi,$ using \eqref{i}, and \eqref{h} we obtain
\begin{align}p^3\overline{\psi}(-4){_4F_3}\left(\begin{matrix}
\overline{\psi},&\phi ,&\phi,&\phi\\
& \phi\overline{\psi},& \phi\overline{\psi},& \phi\overline{\psi}\\
\end{matrix}
\middle|1
\right)&=\frac{S(4,\psi)-2pS(2,\psi)}{g(\psi)^2}\label{j}.
\end{align}
Using \eqref{eqn-Liu}, and \cite[Eq (2.12), Eq (2.13)]{Greene}, we obtain 
\begin{equation}\label{4000}
\sum\limits_{\psi\neq\varepsilon,\phi}\frac{S(2,\psi)}{g(\psi)^2}\binom{\psi}{\phi}^2=p\phi(-1)\sum\limits_{\psi}\binom{\psi}{\phi}^3\psi(1/4)+\frac{1+\phi(-1)}{p^2}.
\end{equation}
Now, substituting $\psi$ by $\phi\psi$ in \eqref{4000} and then using \cite[Eq (2.6)]{Greene}, and Definition \ref{def-1}, we obtain
\begin{equation}\label{4001}
    \sum\limits_{\psi\neq\varepsilon,\phi}\frac{S(2,\psi)}{g(\psi)^2}\binom{\psi}{\phi}^2=(p-1)\phi(-1){_3F_2}\left(\begin{matrix}
        \phi,&\phi, & \phi\\
        &\varepsilon, & \varepsilon
    \end{matrix}
    \middle|\frac{1}{4}
    \right)+\frac{1+\phi(-1)}{p^2}.
\end{equation}
Also, note that for $\psi\neq\varepsilon,\phi$, using \eqref{Jacobi-1}, and \eqref{Jacobi-2}, we have
\begin{equation}\label{4003}
    \frac{1}{g(\psi)^2}\binom{\psi}{\phi}=\frac{\phi(-1)}{p\cdot g(\phi\psi)^2}.
\end{equation}
Finally, substituting \eqref{j} into \eqref{k}, then using \eqref{4001}, and \eqref{4003},  we prove the result. 
\end{proof}

\begin{proof}[Proof of Theorem \ref{corollary-1}]
Putting $k=5$ in Proposition \ref{theorem-2} and substituting the expression from \eqref{s(5)}, we obtain the desired result by comparing the resulting identity with Theorem \ref{theorem-1}.
\end{proof}

\begin{proof}[Proof of Theorem \ref{theorem-5}] 
Using \eqref{eqn-Liu},  and \cite[Eq (2.12), Eq (2.13) ]{Greene}, we obtain 
\begin{equation}\label{4010}
\sum\limits_{\psi\neq\varepsilon,\phi}\frac{S(2,\psi)}{g(\psi)^2}=p\phi(-1)\sum\limits_{\psi}\binom{\psi}{\phi}\psi(1/4)+1+\phi(-1).
\end{equation}
Now using \eqref{Jacobi-1}, \eqref{jacobi-3}, and \eqref{ortho-1} in \eqref{4010}, we obtain
\begin{equation}\label{4011}
    \sum\limits_{\psi\neq\varepsilon,\phi}\frac{S(2,\psi)}{g(\psi)^2}=(p-1)\phi(-3)+1+\phi(-1).
\end{equation}

Finally, substituting \eqref{j}, in Proposition \ref{prop-1}, then substituting \eqref{4011} on the resulting identity, prove our result.
\end{proof}

\begin{proof}[Proof of Theorem \ref{theorem-3}]
For $a\in\F_p^{\times},$ it is easy to verify that $$\sum\limits_{i=0}^{5}\alpha(a)^{5-i}\beta(a)^i=-K(a)^5+4pK(a)^3-3p^2K(a).$$
Using the definition of $T(5,\phi)$ given in \eqref{T(n)} and then using \eqref{Gauss} we have 
\begin{equation}\label{T(5,phi)}
T(5,\phi)=-S(5,\phi)+4pS(3,\phi)-3p^3\phi(-1).
\end{equation}
Now, recall that $N_3(p):=\{(x_1,x_2,x_3)\in({\F}_p^{\times})^3:
\sum\limits_{i=1}^{3}(x_i+x_i^{-1})=0\}.$
Using similar steps as we applied in Proposition \ref{proposition-1}, we obtain
\begin{align}\label{4040}
N_3(p)=&\sum\limits_{a,b\in\F_p}\{ \phi(a^2-1)+1\}\{ \phi(b^2-1)+1\}\{ \phi((a+b)^2-1)+1\}.
\end{align}
Recall that $F(a):=\sum\limits_{x}\phi(x^2-1)\phi((x+a)^2-1)$
as defined in Proposition \ref{proposition-1}. Now, using  \eqref{2.18}, Lemma \ref{lem-2.1}, and  \eqref{id-1}, in \eqref{4040}, thus we obtain 
\begin{equation}
    N_3(p)=\sum\limits_{a}\phi(a^2-1)F(a)+p^2-3p+3.
\end{equation}
Finally, using Lemma \ref{lem-2}, Lemma \ref{lem-2.1}, and Lemma \ref{lem-4} in above expression, we obtain
\begin{equation}\label{99}
    N_3(p)=p(\phi(-3)+\phi(-1))+p^2\cdot{_3F_2}\left( \begin{matrix}
\phi,&\phi,&\phi \\
&\varepsilon,&\varepsilon
\end{matrix}
\middle| \frac{1}{4}\right)+p^2-3p+5.
\end{equation}
On the other hand, if we put $k=3$ in Proposition \ref{theorem-2}  and then substitute the value of $S(3)$ as given in \eqref{S(3)}, we derive
\begin{equation}\label{100}
     N_3(p)=\frac{(p-1)^3}{p}+\frac{(\phi(-3)p^2+2p+1)}{p}+\frac{1}{p}S(3,\phi).
\end{equation}
Comparing \eqref{99}, and \eqref{100}, we obtain 
\begin{equation}\label{S(3,phi)}
   S(3,\phi)= p^2\phi(-1)+p^3\cdot{_3F_2}\left( \begin{matrix}
\phi,&\phi,&\phi \\
&\varepsilon,&\varepsilon
\end{matrix}
\middle| \frac{1}{4}\right).
\end{equation}
Now, substituting \eqref{S(3,phi)} into \eqref{T(5,phi)}, we obtain
\begin{equation}\label{101}T(5,\phi)=-S(5,\phi)+p^3\phi(-1)+4p^4\cdot{_3F_2}\left( \begin{matrix}
\phi,&\phi,&\phi \\
&\varepsilon,&\varepsilon
\end{matrix}
\middle| \frac{1}{4}\right).\end{equation}

For $p\geq7,$ comparing the result of \cite[Theorem 1.2.2,(iii)]{Lim} with \eqref{101}, we have
\begin{equation}\label{102}
    S(5,\phi)=p^3\phi(-1)-p^3\phi(15)+p\phi(-1)a_f(p)^2+4p^4\cdot{_3F_2}\left( \begin{matrix}
\phi,&\phi,&\phi \\
&\varepsilon,&\varepsilon
\end{matrix}
\middle| \frac{1}{4}\right).
\end{equation}

Observe that for a prime $p>3$, we have   $\leg{p}{3}=\phi(-3).$ Finally, if we put $k=5$ in Proposition \ref{theorem-2}, then substituting \eqref{s(5)}, and \eqref{102}, on the resulting identity prove our result.
    
\end{proof}

\begin{proof}[Proof of Theorem \ref{theorem-4}]
It is easy to verify that \begin{equation}\label{110}
    T(6,\phi)=S(6,\phi)-5pS(4,\phi)+6p^2S(2,\phi).
\end{equation}
Using \cite[Theorem 1.1 and Eq. (3.5)]{Dummit}, and \eqref{eqn-Liu}, we have \begin{equation}\label{111}
    S(4,\phi)=-p\cdot a_h(p)-3p^2.
\end{equation}
Substituting \eqref{111} into \eqref{110} and then again substituting the value of $S(2,\phi)$ (using \eqref{eqn-Liu}), we obtain
\begin{equation}\label{112}
    T(6,\phi)=S(6,\phi)+5p^2\cdot a_h(p)+9p^3.
\end{equation}
For $p>7,$ if we compare the result of \cite[Theorem 1.2.2,(iv)]{Lim} with \eqref{112}, then finally we obtain,
\begin{equation}\label{113}
    S(6,\phi)= -a_{\tilde{f}}(p)\cdot p-5p^2\cdot a_h(p)-10p^3.
\end{equation}
Furthermore, if we put $k=6$ in Proposition \ref{theorem-2}, we deduce that 
\begin{equation}\label{N_6(p)}
    N_6(p)=\frac{(p-1)^6}{p}+\frac{1}{p}S(6)+\frac{1}{p}S(6,\phi).
\end{equation}

Now substituting \eqref{S(6)}, and \eqref{113} into \eqref{N_6(p)}, we prove the result.
    
\end{proof}

\begin{proof}[Proof of Corollary \ref{corollary-2}] For $p\geq7,$ comparing Theorem \eqref{theorem-3}, and Theorem \ref{theorem-1}, we finally obtain the identity.
\end{proof}

\section{Acknowledgements}
The first author acknowledges financial support from the
Council of Scientific and Industrial Research (CSIR), Government of India, through a CSIR
Fellowship. The second author is supported by Science and Engineering Research Board/Anusandhan National Research Foundation [CRG/2023/003037].

\subsection*{Data Availability Statement}
No datasets were generated or analyzed during the current study.

\subsection*{Conflict of interest} The author declares that there is no conflict of interest.

	\end{document}